\documentclass[a4paper,11pt]{article}
\usepackage{amsfonts}
\usepackage{amssymb}

\def\F{{\cal F}}

\def\P{{\cal P}}
\def\R{{\cal R}}
\def\S{ {\cal S}}

\def\l{{\lambda}}
\def\s{{\sigma}}
\def\t{{\xi}}

\def\CC{{\mathbb C}}
\def\HH{{\mathbb H}}
\def\RR{{\mathbb R}}
\def\QQ{{\mathbb Q}}

\def\bch{\partial {\cal C}}
\def\dd{ {\partial}}
\def\ML{{\cal  ML}}
\def\MLQ{{\cal  ML}_{\QQ}}
\def\PML{{ \P}{\cal {ML}}}
\def\QF{{\cal Q}{\cal F}}
\def\Sf{ S}
\def\tr{\mathop{\rm Tr}}
 
\newtheorem{thm}{Theorem}[section]
\newtheorem{lemma}[thm]{Lemma}

\newtheorem{cor}[thm]{Corollary}
\newtheorem{prop}[thm]{Proposition}

\newenvironment{proof}{{\sc Proof.}}{$\;\square$ \vskip .2in}

\begin{document}

\title{Thurston's bending measure conjecture  for once punctured torus
groups}
\author{Caroline
Series\\Mathematics Institute\\ Warwick University\\Coventry CV4
7AL,
U.K.}
\date{ }

\maketitle

\bibliographystyle{Plain}

\begin{abstract}   We prove Thurston's bending measure conjecture for
quasifuchsian once punctured torus groups. The conjecture states that
the bending measures of the two components of the convex hull boundary
uniquely determine the group.
\end{abstract}


\section{Introduction}
\label{sec:introduction}

 Thurston conjectured that a convex  
hyperbolic structure on a $3$-manifold with boundary  is uniquely determined by the bending measure on the boundary. In this paper we prove the conjecture in one of the simplest possible examples, namely when the manifold is an interval bundle over a once punctured torus.

In~\cite{KSQF}, the author  and  L. Keen studied the convex hull boundary of this class of hyperbolic $3$-manifolds in great detail, without however addressing Thurston's
conjecture {\em per se}. Since the question came up several times during the Newton Institute programme, it seemed worthwhile to investigate, even though the once punctured torus is a very special case. Given the current state of knowledge, the basic idea of the present proof is rather simple. It does however build on a large number of rather deep results from~\cite{KSQF} and elsewhere.

Let $G \subset PSL(2,\CC)$ be quasifuchsian, so that 
the corresponding quotient $3$-manifold $\HH^3/G$ is an interval bundle over a surface $\Sf$.  The boundary of the convex core of $\HH^3/G$
has two components $\bch^{\pm}$ each of which are pleated surfaces bent along  measured laminations $\beta^{\pm}$ on $\Sf$ called the 
{\em bending measures}. The claim of the bending measure conjecture is  that $\beta^{\pm}$ uniquely determine $G$, up to conjugation in $PSL(2,\CC)$. If the underlying supports of $\beta^{\pm}$ are closed curves, then a result of Bonahon and Otal~\cite{BO} (which applies to general Kleinian groups not just the quasifuchsian case under discussion here) asserts this is indeed the case. They also gave necessary and sufficient conditions for existence, but not uniqueness, of quasifuchsian groups for which the bending measures $\beta^{\pm}$ are any given pair of bending measures $\mu,\nu$. If $\Sf$ is a once punctured torus, these conditions reduce  to $i(\mu,\nu) > 0$ and the proviso that the weight of any closed curve is less than $\pi$. In a recent preprint~\cite{Bon}, Bonahon has also shown uniqueness for general  quasifuchsian groups which are sufficiently close to being Fuchsian.

The object of~\cite{KSQF} was to describe the space
$\QF$ of quasifuchsian groups when the surface $\Sf$ is a once punctured torus in terms of the geometry of $ \bch^{\pm}$. This was done  by analysing the
{\em pleating planes} $\P(\mu,\nu)$ consisting of all groups 
whose bending laminations lie in a particular pair of projective
classes $[\mu],[\nu]$, see Section~\ref{sec:prelims}.  In place of the bending measures, we concentrated on the lengths $l_{\mu}$ and $l_{\nu}$ of $\mu$ and $\nu$, which can be extended to well-defined  holomorphic functions, the {\em complex lengths} $\l_{\mu}$ and $\l_{\nu}$, on $\QF$, see Section~\ref{sec:plvarieties}. We showed, see Theorem~\ref{thm:plplanes}, that $\l_{\mu}$ and $\l_{\nu}$
are local holomorphic coordinates for $\QF$ in a neighbourhood of $\P(\mu,\nu)$
whose restrictions to $\P(\mu,\nu)$ are real valued. In this way we
obtained a diffeomorphism of $\P(\mu,\nu)$ with a certain open 
subset in $\RR^+ \times \RR^+$. In particular,
the lengths $l_{\mu}$ and $l_{\nu}$
 uniquely determine the group.

Given this background, to prove the bending measure  conjecture
we just need to show that, restricted to a given pleating variety, the map from lengths to bending measures is injective. (Here by bending measure, we really mean the scale factors which relate  $\beta^{\pm}$ to a fixed choice of laminations $\mu \in [\mu], \nu \in [\nu]$.) For rational laminations, as mentioned above, we already know by~\cite{BO} that angles determine the group. In~\cite{ChS}, again in the context of  general Kleinian groups, we showed that for rational bending laminations the map from
bending angles to lengths is injective and moreover that its Jacobian is  symmetric and negative definite. If $\Sf$ is a once punctured torus, it follows (see Section~\ref{sec:monotone}) that if both bending laminations are rational and if we fix the length
of the bending line on say $\bch^{+}$, then the bending angle on $\bch^{-}$
is a monotonic function of the length on $\bch^{-}$.

Now keeping the bending line on $\bch^{+}$ fixed and of fixed length, we 
take limits as the the bending laminations on $\bch^{-}$ converge projectively to an arbitrary irrational lamination $\nu$. Using the fact that the limit of monotone functions is monotone, we deduce the scale factor of the bending measure on $\nu$ is still a monotone function of the length $l_{\nu}$. Since this scale factor is  real analytic, it is either strictly monotonic or constant; a global geometrical argument rules out constancy. An elaboration of the argument in Sections~\ref{sec:constantangle} and~\ref{sec:mainproof} then allows us relax the requirement that the bending lamination and length  on $\bch^{+}$ remain fixed and prove the conjecture.

In~\cite{KSQF} we also showed that the rational pleating varieties for which the supports of the bending laminations are simple closed curves are dense. This is not enough for the present proof: to compare our monotone functions properly, we need the ``limit pleating theorem'' of~\cite{KSQF}, see Theorem~\ref{thm:limitpleating}. This is a deep result closely related to the ``lemme de fermeture'' in~\cite{BO}. Roughly, it asserts the existence of an algebraic limit in $\P(\mu,\nu)$ for any sequence along which the pleating lengths $l_{\mu},l_{\nu}$ remain bounded.

A simpler version of the same analysis, not spelled out here, would show  that the bending angle also uniquely determines the group for the so-called Riley slice of Schottky space, see~\cite{KSRiley}.

We would like to thank in particular Cyril Lecuire and Pete Storm whose
questioning fixed this problem in my mind.

\section{Preliminaries}
\label{sec:prelims}

Throughout the paper,  $\Sf$ will denote a fixed topological once punctured torus.
A representation $\pi_1(\Sf) \to  PSL(2,\CC)$ is called quasifuchsian if
the image is discrete and torsion free, if the image of a simple loop around the puncture is parabolic, and if the quotient hyperbolic manifold $M=\HH^3/G$ is homeomorphic to $\Sf \times (-1,1)$. It is Fuchsian if it is quasifuchsian and if in addition the representation is conjugate to a representation into $PSL(2,\RR)$.
Let $\F$ and $\QF$ denote  the spaces of Fuchsian and quasifuchsian representations respectively, modulo conjugation
in $PSL(2,\CC)$.
The space $\QF$ is a smooth complex manifold of dimension $2$, with   natural holomorphic structure
induced from $PSL(2,\CC)$.

\subsubsection{Measured laminations} We assume the reader is  familiar  with geodesic laminations, see for example~\cite{ThuN, OtalH}.
A  measured lamination  $\mu$ on $\Sf$ consists of a geodesic
lamination,
called the support of $\mu$ and denoted $|\mu|$, together
with a transverse measure, also denoted $\mu$.  
We topologise the space $\ML$ of
all measured laminations on $\Sf$ with the topology of weak convergence of transverse measures, that is, two laminations are close if the measures they assign to any finite set of transversals are close. 
We write $[\mu]$ for the projective class of $\mu
\in \ML(\Sf)$ and denote the set of projective equivalence classes on $\Sf$ by $\PML$.

Let $\cal S$ be the set of simple closed curves on $\Sf$.
We call $\mu \in \ML$ {\em
rational} if $|\mu| \in \S$. (On a more general surface, we say a lamination is rational if its support is a disjoint union of simple closed curves.) Equivalently, $\mu$ is rational if $\mu = c\delta_{\gamma}$ where $\gamma \in \S, \ c \ge 0$ and  $\delta_{\gamma}$  is the measured lamination   which assigns unit mass to each intersection
with $\gamma$.  The
set of all rational measured laminations on $\Sf$ is denoted
$\MLQ$; this set is dense in $\ML$.

The geometric intersection number $i(\gamma,\gamma')$ of two
geodesics
$\gamma,\gamma' \in {\cal S}$ extends to a jointly continuous function
$i(\mu,\nu)$
on $\ML$. 
It is special to the once punctured torus that
$i(\mu,\nu) > 0 $ is equivalent to $[\mu] \neq [\nu]$, moreover (since all laminations on $\Sf$ are uniquely ergodic) if $\mu,\nu \in \ML$ with $|\mu|=|\nu|$, then $[\mu]=[\nu]$.

For general surfaces, convergence of laminations in the topology of $\ML$ does not imply Hausdorff convergence of their supports. However:

\begin{lemma}(\cite{KSQF} Lemma 1) 
\label{lemma:hausdorff}
Suppose $\Sf$ is a once punctured torus. Suppose that $\mu \notin \MLQ$ and that $\mu_n \to \mu$ in the topology of weak convergence on $\ML$. Then $|\mu_n| \to |\mu|$
in the Hausdorff topology on the set of closed subsets of $\Sf$.
\end{lemma}

\subsubsection{Lengths}

Given a hyperbolic structure on $\Sf$ (associated to a Fuchsian representation of $\pi_1(\Sf)$) and  $\gamma \in \pi_1(\Sf)$, we define the length $l_{\gamma}$ to  be the hyperbolic length of the unique geodesic  freely homotopic to $\gamma$. This definition can be extended to general laminations $\mu \in \ML $. The following theorem summarizes the results of~\cite{KerckEA}, Lemma~2.4 and \cite{KerckN}, Theorem~1:
\begin{prop}
\label{prop:Kerck2}
The function $( c\delta_{\gamma},p) \mapsto c
l_{\delta_{\gamma}}(p)$ from
$\ML_{\QQ} \times \F$ to $\RR^+$ extends to a real analytic 
function
$(\mu,p)
\mapsto l_{\mu}(p)$ from $\ML \times \F $ to $ \RR^+$. 
 If $\mu_n
\in \ML_{\QQ}$,
$\mu_n \to \mu$ then  $ l_{\mu_n}(p) \to
l_{\mu}(p)$
uniformly on compact subsets of $\F$.
\end{prop}

We showed in~\cite{KSQF} that for
 $\mu \in \ML$, the length function $l_{\mu}$ on $\F$ extends to
a  non-constant holomorphic function $\lambda_{\mu}$, called the {\em complex length} of  $\mu$, on $\QF$. The extension is done in such a way that $\lambda_{c\mu} = c \lambda_{\mu}$ for $c>0$,  and such that if $\mu = \delta_{\gamma} $ then $q \mapsto \lambda_{\mu}(q) $ is a well-defined branch of the complex length $ 2{\cosh}^{-1}  \tr \rho(\gamma)/2$, where 
$\rho: \pi_1(\Sf)\to PLS(2,\CC)$ represents $q \in \QF$.  
\begin{prop}(\cite{KSQF} Theorem 20)
\label{prop:cxlength} 
The family  $\{\lambda_{\mu}\}$ is uniformly bounded and equicontinuous on compact subsets of $\ML \times \QF$, in particular if $\mu_n \to \mu \in \ML$ and $q_n \to q \in \QF$ then 
$ \lambda_{\mu_n}(q_n) \to \lambda_{\mu}(q) $.
 \end{prop}

The real part $l_{\mu}$ of the complex length $\l_{\mu}$ is the lamination length of $\mu$ in $\HH^3/G$ as defined in~\cite{ThuN} p. 9.21.

\subsubsection{The Thurston boundary}     

 We recall the fundamental inequality which governs Thurston's compatification of 
Teichm\"uller space $\F$ with projective measured lamination space $\PML$, see~\cite{FLP} Lemme II.1 and~\cite{ThuII} Theorem 2.2:
\begin{prop}
\label{prop:thurstonbndry} Suppose $\s_n \to [\xi] \in \PML$. Then there exist $C>0$, $d_n \to \infty$ and $\xi_n \to \xi$ in $\ML$ such that
$$d_n i(\xi_n,\zeta) \le l_{\zeta}(\s_n) \le d_n i(\xi_n,\zeta) + C  l_{\zeta}(\s_0) $$ for all $\zeta\in \ML$.
\end{prop}

\begin{cor}
\label{cor:thurstonbndry}  Suppose given a sequence of measured laminations $\mu_n$ on the once punctured torus $\Sf$
such that $\mu_n  \to \mu$ in $\ML$,  and suppose that $\s_n \in \F$ is a sequence of metrics such that $l_{\mu_n}(\s_n) \to 0$. Then $\s_n \to [\mu] \in \PML$.
\end{cor}
\begin{proof} If $\s_n$ converged to a point $\s_{\infty}$ in $\F$ then $l_{\mu}(\s_{\infty}) = \lim_n \  l_{\mu_n}(\s_n) = 0$ which is impossible. Thus $\s_n$ converges to some point $  [\xi] \in \PML$. Substituting in the fundamental inequality we obtain
$\xi_n \to \xi$ in $\ML$ and $d_n \to \infty$
such that:
$$d_n i(\xi_n,\mu_n) \le l_{\mu_n}(\s_n) \le d_n i(\xi_n,\mu_n) + C  l_{\mu_n}(\s_0).$$

Since $i(\xi_n,\mu_n) \to i(\xi,\mu)$, we see that $l_{\mu_n}(\s_n)/d_n \to i(\xi,\mu)$.
We conclude that $i(\xi,\mu) = 0$, and hence, since $\Sf$ is a punctured torus, that $[\xi]=[\mu]$.
\end{proof}

\subsubsection{Bending measures} 
Let $q \in \QF$ and let $G=G(q)$ be a group representing $q$.  The convex core $\cal C$ of $\HH^3/G$ is the smallest closed set containing all closed geodesics; it is the projection to the quotient of the convex hull of the limit set and has non-zero volume if and only if 
$G \in \QF - \F$. In this case its boundary  has two connected components $ \bch^{\pm}$ each of which is homeomorphic to $\Sf$.
The metric induced on $  \bch^{\pm}$ from $\HH^3/G$
makes each a pleated surface.
The  bending laminations of these surfaces 
 carry  natural transverse measures,
 the {\em bending measures} $\beta^{\pm}$, 
see~\cite{EpM, KSconvex}. The underlying geodesic laminations $|\beta^{\pm}|$
 are called the {\em pleating loci} of $G$.
Since the same geodesic lamination cannot be the pleating locus on both sides, we have $i(\beta^+,\beta^-) >0$, see~\cite{KSQF,BO}. 
If $ \beta^+$ is rational with support $\gamma \in \S$, then $ \beta^+ = \theta_{\gamma} \delta_{\gamma}$ where $ \theta_{\gamma}$ is the bending angle along $\gamma$. In this case we have the obvious constraint  $ \theta_{\gamma} < \pi$. If $G$ is Fuchsian we define $\beta^{\pm}= 0$. One of the main results of~\cite{KSconvex} is that $\beta^{\pm}$ are continuous functions on $\QF$. 

The following central existence result is a special case of Bonahon and Otal's ``Lemme de
fermeture'':
\begin{thm}(\cite{BO} Proposition 8)
\label{thm:bonotal} Suppose $\mu,\nu \in \ML$ with $i(\mu,\nu)>0$. If $\mu = c \delta_{\gamma}$ is rational assume also that $c<\pi$,  and similarly for $\nu$. Then there exists a quasifuchsian group $G \in \QF$ such that $\beta^+ =  \mu$ and  $\beta^- = \nu$.
If $\mu,\nu$ are rational, then $G$ is unique.
\end{thm}

The object of this paper is to prove the uniqueness part of this
statement for arbitrary $\mu$ and $\nu$.

\section{Pleating varieties} 
\label{sec:plvarieties}

In this section we review the results we shall need from~\cite{KSQF}.
 
Fix $\mu,\nu \in \ML$.
The {\em $(\mu,\nu)$-pleating variety} is the set $\P(\mu,\nu) \subset \QF $ such that 
$\beta^+ \in [\mu]$ and $\beta^- \in [\nu]$ with $\beta^{\pm} \neq 0$. 
(The last condition is equivalent to $\P(\mu,\nu) \cap \F = \emptyset$.)
Notice that $\P(\mu,\nu)$ actually only depends on 
the projective classes $[\mu]$ and $[\nu]$.  By Theorem~\ref{thm:bonotal} or alternatively~\cite{KSQF} Theorem 2, $\P(\mu,\nu)$ is non-empty if and only if 
 $i(\mu,\nu)>0$.
By~\cite{KSQF} Proposition 22, the complex length
 $\lambda_{\mu} $ is real valued whenever the projective class of $\beta^{+}$ (or $\beta^-$) is $[\mu]$. In this case, $l_{\mu} = \Re \l_{\mu}$
is both the lamination length in $\HH^3/G$ and the length of $\mu$
in the hyperbolic structure  on $\bch^{+}$  (or $\bch^{-}$).

 Let $(\mu,\nu)$ be measured laminations on $\Sf$ with
$i(\mu,\nu)>0$. Then by~\cite{KerckLM}, the function $l_{\mu} + l_{\nu}$ has a unique minimum on $\F$. Moreover as $c$ varies in $(0,\infty)$, the minima of $l_{\mu} + l_{c\nu}$ form a line $\cal L(\mu,\nu)$; it is easy to see that this line only depends on the projective classes $[\mu]$ and $[\nu]$. The lines $\cal L(\mu,\nu)$ vary continuously with
$\mu$ and $\nu$, see~\cite{KerckLM}.
(In~\cite{KSQF} we gave a slightly different description of $\cal L(\mu,\nu)$ in terms of the minimum of $l_{\mu}$ on an earthquake path along $\nu$, namely, $\cal L(\mu,\nu)$ is the locus where $\frac{\dd l_{\mu}} {\dd t_{\nu}}$ vanishes, where $t_{\nu}$ denotes the earthquake along $\nu$. This definition only works for the once punctured torus; for the equivalence of the two definitions, see~\cite{KSQF} Lemma 6.) We proved in~\cite{SmB} Theorem 1.7 (in the context of  arbitrary hyperbolisable surfaces $\Sf$) that the closure of $\P(\mu,\nu)$ meets $\F$ exactly in  $\cal L(\mu,\nu)$.

For each $c>0$, the line $\cal L(\mu,\nu)$ meets the horocycle $l_{\mu} = c$ in $\F$ in exactly one point.  Let $f_{\mu,\nu}: \RR^{+} \to
\RR^{+}$  be the function $f_{\mu,\nu}(c) = l_{\nu} (d)$ where $d = {\cal L}(\mu,\nu) \cap l_{\mu}^{-1}( c)$. We showed (\cite{KSQF} Lemma 6) that $f_{\mu,\nu}$ is surjective and strictly monotone decreasing.
We denote by $\R(\mu,\nu)$  the
open region in $\RR^{+} \times
\RR^{+}$ bounded between the two coordinate axes and 
the graph of $f_{\mu,\nu}$.
The following description of $\P(\mu,\nu)$ is one of the main results of~\cite{KSQF}:
\begin{thm}(\cite{KSQF} Theorem 2) 
\label{thm:plplanes}
Let $\mu,\nu \in \ML$ be measured laminations with
$i(\mu,\nu)>0$. 
Then the function $\l_{\mu} \times \l_{\nu}$ is locally injective in a neighbourhood of 
 $\P(\mu,\nu)$, moreover its restriction to  $\P(\mu,\nu)$ is a 
diffeomorphism to $\R(\mu,\nu)$.
\end{thm}

This means that $\P(\mu,\nu)$ is totally real, in other words, there are local holomorphic coordinates such that 
a neighbourhood of $x \in \P(\mu,\nu) \hookrightarrow \QF$ is identified with a neighbourhood
of $0 \in \RR^2 \hookrightarrow \CC^2$. 
Moreover $\P(\mu,\nu)$ is a connected real $2$-manifold, on which we can take $(l_{\mu} ,l_{\nu})$
as global real analytic coordinates, where as above we define 
$l_{\mu} = \Re \l_{\mu},l_{\nu} = \Re \l_{\nu}$.

\begin{cor}
\label{cor:line}
The set $L_c = l_{\nu}^{-1}(c) \subset \P(\mu,\nu)$ is connected and hence can be regarded as a line in $\R(\mu,\nu)$ parameterised by $l_{\mu}$. 
\end{cor}
\begin{proof}
This is immediate from the fact that $f_{{\mu},{\nu}}$ is strictly monotonic.
\end{proof}

 A crucial ingredient of Theorem~\ref{thm:plplanes} was the following, which we call the {\em limit pleating theorem}:
\begin{thm}(\cite{KSQF} Theorem 5.1)
\label{thm:limitpleating}
Suppose $q_n \in \P(\mu,\nu)$  is a sequence such that the lengths $ l_{\mu}(q_n),l_{\nu}(q_n)$ are
uniformly bounded above. Then there is a subsequence of 
$q_n$ such that the corresponding groups $G(q_n)$ converge algebraically to a group   
$G_{\infty}$. If both $\lim_n l_{\mu}(q_n)$ and $ \lim_n l_{\nu}(q_n) $  are strictly positive,
then $G_{\infty}\in \F \cup \P(\mu,\nu)$ and  the convergence is strong.
\end{thm}

It is convenient to have a somewhat stronger version of this theorem, in
which we allow the bending loci to vary also:
\begin{thm}
\label{thm:limitpleating1}
Suppose that $\mu,\nu \notin \MLQ$ and that 
$\mu_n \to \mu$ and $\nu_n \to \nu$.
Suppose $q_n \in \P(\mu_n,\nu_n)$ and that $l_{\mu_n}, l_{\nu_n}$ are
uniformly bounded above. Then there is a subsequence of 
$q_n$ such that $G(q_n)$  converges algebraically to a group   
$G_{\infty} $. If both $\lim_n l_{\mu_n}(q_n)$ and $ \lim_n l_{\nu_n}(q_n) $  are strictly positive,
then $G_{\infty}\in \F \cup \P(\mu,\nu)$ and the convergence is strong.
\end{thm}

The proof is almost identical to that in~\cite{KSQF}.
We make the assumption that $\mu,\nu \notin \MLQ$ since we need that the convergence of $\mu_n$ to $\mu$ be Hausdorff; this follows from
Lemma~\ref{lemma:hausdorff}.

Very briefly,  the proof of Theorem~\ref{thm:limitpleating} goes as follows. The bounds on $l_{\mu}$ and $l_{\nu}$  allow  one to use Thurston~\cite{ThuII} Theorem 3.3 (as in the proof of the double limit theorem) to conclude the existence of  an algebraic limit $G_{\infty}$. The main part of the work is to show that 
$G_{\infty}$ is quasifuchsian.  One uses continuity of the lamination length to show that the laminations $\mu$ and $\nu$ are realised in the limit $3$-manifold $\HH^3/G_{\infty}$. By carefully using the full force of the algebraic convergence, one can then show that, in the universal cover $\HH^3$,  the lifts of the geodesic laminations $|\mu|$ and $|\nu|$ from $\HH^3/G(q_n)$  approach  the lifts of the realisations of $|\mu|$ and $|\nu|$ from $\HH^3/G_{\infty}$ in the Hausdorff topology. 
From this it is not hard to deduce that the laminations $|\mu|$ and $|\nu|$ lie in boundary of the convex core of $\HH^3/G_{\infty}$. One deduces that this boundary has two components (or one two-sided component if $G_{\infty}$ is Fuchsian) from which it follows that $G_{\infty} \in \QF$.
Inspection shows that exactly the same proof gives the more general version Theorem~\ref{thm:limitpleating1},
provided we have Hausdorff convergence of $|\mu_n|$ to $|\mu|$ and 
$|\nu_n|$ to $|\nu|$.

Alternatively, one can modify the proof of the Lemme de fermeture~\cite{BO}  Proposition 8. They make an assumption on the limit bending laminations (which in particular rules out that the limit group is Fuchsian) but since  the main use of this assumption is to get a  bound  on the lengths $l_{\beta^{\pm}}$, their proof can be adapted relatively easily to our case.

\subsubsection{Scaling functions}

For fixed $\mu,\nu \in \ML$ we define the {\em scaling functions} $\t_{\mu}, \t_{\nu}$ by  $\beta^+ = \t_{\mu} \mu$ and $\beta^- = \t_{\nu} \nu$.

\begin{prop}
The scaling function  $\t_{\mu}$ is real analytic on 
$\P(\mu,\nu)$.  
\end{prop}
\begin{proof} Suppose first that $\mu \in \MLQ$ and 
let $\gamma = |\mu|$. Let 
$\lambda_{\gamma}, \tau_{\gamma }$ be complex Fenchel Nielsen coordinates for $\QF$ with respect to 
$\gamma$. Here $\lambda_{\gamma}$ is the complex length of $\gamma$ and $\tau_{\gamma }$ is the complex twist 
along $\gamma$, see for example~\cite{KSQF} or~\cite{KSbend}. The quasifuchsian group $G$ can be described up to conjugation by these parameters, which are holomorphic functions on $\QF$. As explained in detail in~\cite{KSQF},
if $\beta^+ \in [\mu]$ then  $\beta^+ = \theta_{\gamma} \delta_{\gamma}$ and  $\theta_{\gamma} = \Im \tau_{\gamma }$. This easily gives the result.

For general $\mu$, we can obtain any group for which $\beta^+ \in [\mu]$ by a complex earthquake along $\mu$. The scale factor $\t_{\mu}$ is easily seen to be  
 $\Im \tau_{\mu}$. This function is again holomorphic,
see Section 7.3 of~\cite{KSQF}. One can also obtain this result by using the shear-bend coordinates developed by Bonahon in~\cite{BonS}, see also~\cite{Bon}.
\end{proof}

In the special case under discussion the bending measure conjecture can thus be viewed as the assertion that the map $\P(\mu,\nu) \to \RR^+ \times \RR^+$,  $ (l_{\mu},l_{\nu}) \mapsto (\t_{\mu},\t_{\nu}) $ is injective.

\section{Global geometry}
\label{sec:geometry}

We shall need various results which control the behaviour of bending angles versus bending lengths.
The following basic inequality is due to Bridgeman: 
\begin{prop}(\cite{Br} Proposition 2)
\label{prop:bridgeman}
There exists a universal constant $K>0$ such that if $\beta^+ \in [\mu]  \in \PML$, then 
$l_{\mu}\t_{\mu} \le K$.
\end{prop}

Say $G \in \QF$.
If $\gamma \in \pi_1(\Sf)$ we let $\gamma^*$
be the geodesic representative  of $\gamma$ in $\HH^3/G$
and we write $\gamma^{\pm}$ for the geodesic representative in the hyperbolic structure of $\bch^{\pm}$. Denote the respective lengths by 
$l_{\gamma^*}$ and $l_{\gamma^{\pm}}$. We recall that $l_{\gamma^*} = \Re \l_{\gamma^*}$ where $\l_{\gamma^*}$ is the complex length defined in Section~\ref{sec:plvarieties}.
The following result, a simplified version of Lemma A.1 of~\cite{Lec}, see also Proposition 5.1 in~\cite{SmB}, gives a good comparison between the lengths of $\gamma^+$ and $\gamma^*$ 
when the intersection of $\gamma$ with the bending measure is reasonably small. 
\begin{prop}
\label{prop:lecuire} Let $\beta^+$ be the bending measure of the component $\bch^{+}/G $ of the convex hull boundary of 
the manifold $\HH^3/G$. Let $\gamma$ be a simple closed curve with geodesic representatives $\gamma^*$ in $\HH^3/G$ and $\gamma^+$ in $\bch^{+}$. Then there exist  $A,B>0$ such that if $i(\gamma,\beta^+)  < \pi/12$, then
$ l_{\gamma^+} <  Al_{\gamma^*} +B$.
\end{prop}

Now we can prove a result which will give the global control we need.
\begin{prop}
\label{prop:lengthbounds}
Fix $a>0$. Then for any $K>0$, there exists $\epsilon>0$
such that whenever $q  \in \P(\mu,\nu)$ and $l_{\mu}(q)< \epsilon $ with
$\t_{\mu} (q)= a $,   then 
$l_{\nu} (q)> K $.
\end{prop}
\begin{proof}  
If the result is false, then there is a sequence 
$q_m \in \P(\mu,\nu)$ for which $\t_{\mu} (q_m)= a $
and  $l_{\nu} (q_m) \le K $ but 
$l_{\mu}(q_m) \to 0$.
By  Theorem~\ref{thm:limitpleating}, under these conditions the groups $G(q_m)$ have an algebraic limit so that in particular, the hyperbolic lengths 
$l_{\gamma^*}(q_m) = \Re \l_{\gamma^*}(q_m)$ are uniformly bounded above for any simple
closed curve $\gamma \in \pi_1(\Sf)$.

Suppose first that $\mu \in \MLQ$, $\mu = c \delta_{\eta}$ for $\eta \in \S$. Then the bending angle $\theta_{\eta}$ is $c\t_{\mu}$, where necessarily $\theta_{\eta} < \pi$. Since $\Sf$ is a punctured torus, we can choose a simple closed curve $\gamma$ such that $i(\gamma, \eta) =1$. As on p.195 in~\cite{PS}, for any $G \in \QF$, we have the equation
$$\cosh \tau_{\eta}/2 = \cosh \l_{\gamma^*}/2 \; \tanh \l_{\eta^*}/2$$ where $\tau_{\eta}$ is the complex Fenchel Nielsen twist along $\eta$.
Since $\eta$ is assumed to be the bending locus of $\bch^{+}$, we have $\tau_{\eta} = t_{\eta} + i\theta_{\eta} $ (see~\cite{KSQF} or~\cite{KSbend} for detailed discussion). 
Noting that $\l_{\eta^*} = l_{\eta^*} \in \RR$ and taking real parts we see that  
 if  $\theta_{\eta}=ac < \pi$ is fixed, then $l_{\eta^*} \to 0$ implies that  $l_{\gamma^*} = \Re \l_{\gamma^*} \to \infty$. This contradicts the uniform upper bound to the lengths $l_{\gamma^*}(q_m)$ from the first paragraph.

Now suppose that $\mu \notin \MLQ$ and that $\t_{\mu} =a$ is fixed. By following a a leaf of $|\mu|$ which returns sufficiently close to itself along some transversal, we can choose a simple closed curve $\gamma$ such that $  i(\gamma, \beta^+) = ai(\gamma, \mu) < \pi/12$.  
Since $\mu$ is the bending measure of $\bch^{+} (q_m)$, the lengths $l_{\mu^*}$ and $l_{\mu^+}$ of 
the lamination ${\mu}$ in $\HH^3/G(q_m)$ and on the pleated surface $\bch^{+} (q_m)$
coincide. 
Since $l_{\mu^+} (q_m) \to 0$, by Corollary~\ref{cor:thurstonbndry}
the hyperbolic structures of the pleated surfaces $\bch^{+} (q_m)$
tend to $[\mu] \in \PML$. 
Since $i(\gamma, \mu)>0$, it follows that  $l_{\gamma^+} (q_m) \to \infty$. By Proposition~\ref{prop:lecuire}, we have 
$l_{\gamma^*}(q_m) \to \infty$, and the same  contradiction as before completes the proof.
\end{proof}

\section{Monotonicity of angle for fixed length}
\label{sec:monotone}

Suppose that the bending laminations $\mu,\nu$ are rational,
supported by simple closed curves $\gamma,\delta \in \S$.
Theorem~\ref{thm:plplanes} asserts that the pleating variety $\P(\mu,\nu) = \P(\gamma,\delta) $ is an open  real $2$-manifold parameterised by the lengths $l_{\gamma}, l_{\delta}$.
 Theorem~\ref{thm:bonotal} 
shows on the other hand   that
  the bending angles $\theta_{\gamma}, \theta_{\delta}$ are global coordinates for $\P(\gamma,\delta)$.

In~\cite{ChS} we studied the relationship between bending angle and length for general hyperbolic $3$-manifolds with boundary, where the bending laminations were rational.
In particular, we showed in Proposition 7.1, that if we regard the lengths $l_i$ of the bending lines as functions of the bending angles $\theta_i$, then the Jacobian matrix 
$  \left( \frac{\partial l_{i}}{\partial 
\theta_{j}} \right)$  evaluated at any point in the relevant pleating variety
is negative definite and symmetric.
(The main point of~\cite{ChS} is to establish that in general lengths are parameters. Given  that both lengths and angles are parameters, so that  the Jacobian  is non-singular, the fact that it is symmetric and negative definite follows
 easily  from the  Schl\"afli formula for variation of volume of the convex core and using the symmetry of second derivatives, see~\cite{ChS}.) 

Since the inverse of a negative definite symmetric matrix is also negative definite and symmetric, and since
the diagonal entries of a negative definite matrix cannnot vanish,
we deduce in our special case that $\frac{\partial \theta_{\gamma}}{\partial 
l_{\gamma}}(q)<0$ for $q \in \P(\gamma,\delta)$ (where the partial derivative is taken keeping $l_{\delta}$ fixed). Recall from Corollary~\ref{cor:line} that $L_c = l_{\nu}^{-1}(c) \subset \P(\mu,\nu)$  can be regarded as a line parameterised by $l_{\mu}$. 
\begin{prop}\label{prop:ratmonotone}
 Let $\mu, \nu \in \ML_{\QQ}$. Fix
$c>0$. Then the scaled bending angle $\t_{\mu}$ is a strictly monotone function of $l_{\mu}$ on the line  $L_c$. 
\end{prop}

We want to take limits to prove the same result for general $\mu, \nu$.
We use the following simple fact about analytic functions:
\begin{lemma} Let $f_n$ be a sequence of real valued monotonic functions
on $(a,b) \subset \RR$ which converges pointwise to 
a real analytic function $f:(a,b) \to \RR$. The $f$ is either
strictly monotonic or constant.
\end{lemma}

Notice that  a {\em real analytic} function  $f$ may be strictly monotonic even though its derivative vanishes at some points. (The sequence $f_n(x) = x^3 + x/n$ is an instructive example.)  
Thus in the following result we only claim monotonicity and not necessarily that
$\frac{\partial \t_{\mu}}{\partial 
l_{\mu}}<0$.

\begin{prop}\label{prop:linemonotone} Let $\mu, \nu \in \ML$. Fix
$c>0$. Let 
 $L_c = l_{\nu}^{-1}(c) \subset \P(\mu,\nu)$.
 Then the function $\t_{\mu}$ is a strictly monotone function of
$l_{\mu}$ on the line  $L_c$. 
 \end{prop}
 \begin{proof}
Choose $\mu_n, \nu_n \in \MLQ$ with $\mu_n \to \mu$ and $\nu_n \to \nu$
in $\ML$. If either $\mu$ or $\nu$ is in $\MLQ$ then choose
the corresponding sequence to be constant. By Lemma~\ref{lemma:hausdorff}, both  sequences
of laminations also converge in the Hausdorff topology.

Fix $b,c>0$ such that $ (b,c) \in \R(\mu,\nu)$. 
 By Theorem~\ref{thm:plplanes}, there is a unique point $w (b,c)$ with 
$l_{\mu} = b$ and $l_{\nu} = c$ in $\P(\mu,\nu)$.
Since the lines of minima vary continuously with the laminations,
for  sufficently large $n$,  $\R(\mu_n,\nu_n)$
is close to $\R(\mu,\nu)$ and hence there is a unique point
$w_n(b,c)$ in  $\P(\mu_n,\nu_n)$ for which $l_{\mu_n} = b$ and $l_{\nu_n} = c$. 
 By the strengthened version of the limit pleating theorem~\ref{thm:limitpleating1}, up to extracting a subsequence these points converge to a point in $ \P(\mu,\nu)$. At this point
 $l_{\mu} = b$ and $l_{\nu} = c$, hence by Theorem~\ref{thm:plplanes} the limit point is unique and must equal $w (b,c)$. Now keep $c$ fixed and vary $b$. For each $n$, we know from Proposition~\ref{prop:ratmonotone}  that $\t_{\mu_n}(w_n(b,c)) $ is a
strictly monotone decreasing function of $b$.  By continuity we have
$\t_{\mu_n}(w_n(b,c)) \to  \t_{\mu}(w (b,c))$. Hence for $c$ fixed, 
$\t_{\mu}(w (b,c))$ is monotone decreasing.
Moreover for fixed $l_{\nu }= c$, the function $\t_{\mu}$ is real analytic
in  $l_{\mu}$. Hence it is either strictly decreasing or constant.
However $\t_{\mu}$ is certainly not constant since 
$\t_{\mu} \to 0 $ as we approach Fuchsian space along 
$L_c$.
The result follows.
\end{proof}

We remark that a somewhat simpler version of the above proof would show  monotonicity of angle in the Riley slice of Schottky space, proving the bending measure conjecture for a genus two handlebody with both handles pinched (so that the boundary is a sphere with four punctures).

\section{The constant angle variety}
\label{sec:constantangle}

Fix $a>0$ and let $V_a= \{q \in \P(\mu,\nu) | \t_{\mu} =a \}$.
Notice that $V_a \neq \emptyset$ by Theorem~\ref{thm:bonotal}.

\begin{lemma}
\label{lem:notisolated} $V_a$ contains no isolated points. 
\end{lemma}
\begin{proof} Suppose the contrary, and let $q \in \P(\mu,\nu)$
be an isolated point of $V_a$. Then there exists an open disk $D$
containing $q$ with $V_a \cap D = \{q\}$. The image of $D-\{q\}$ under
$\t_{\mu}$ is connected, and it follows that $a$ is a local maximum
or minimum for the function $\t_{\mu}$ restricted to the line $L_c =
l_{\nu}^{-1}(c)$, $c = l_{\nu}(q)$. This contradicts the strict
monotonicity of $\t_{\mu}$ on $L_c$ from Proposition~\ref{prop:linemonotone}.
\end{proof}

\begin{lemma}
\label{lem:dim1} $V_a$ is a real algebraic variety of dimension
$1$. 
\end{lemma}
\begin{proof} By definition, $V_a$ has dimension $2-r$ where 
$r$ is the maximum rank of $(\frac{\dd
\t_{\mu}}{\dd l_{\mu}},\frac{\dd \t_{\mu}}{\dd l_{\nu}} )$ on
$V_a$. Clearly $r \le 1$. If $r = 0$  then the Jacobian vanishes
identically on $V_a$, from which we deduce that all partial derivatives
of $\t_{\mu}$ vanish identically. This would mean  the real analytic function
$\t_{\mu}$ was constant in a neighbourhood of $V$ in $\P(\mu,\nu)$, contradicting Proposition~\ref{prop:linemonotone}.
\end{proof}

We can now use the classical local description of real analytic
varieties, see for example Milnor~\cite{Mi} Lemma 3.3: 
\begin{thm}
\label{thm:milnor}   Let $P = (x_0,y_0)$ be a non-isolated point in a real
one dimensional algebraic variety $V \subset \RR^2$. Then there is  a
neighbourhood 
of $P$ in which $V$ consists of a finite number of branches, each of
which is homeomorphic to an open interval on $\RR$.
After interchanging the two coordinates if necessary, 
we may assume this homeomorphism has the form
$ t \mapsto (x_0+   t^k, y_0+ \sum_{r=1}^{\infty} b_r t^r )$ where  the highest common factor of $k$ and the indices of the non-vanishing coefficients $b_r$ is $1$.  
\end{thm}

\begin{prop}
\label{prop:shapeofV}
 $V_a$ consists of a unique connected component with no branch points
which
intersects each line $L_c$ in a unique point.
\end{prop}

\begin{proof} We already know by Proposition~\ref{prop:linemonotone}
that $V_a \cap L_c$ consists of at most one point. In other words, the
restriction of $l_{\nu}$ to $V_a$ is injective, hence in particular, has
no local maxima or minima. Denote this restriction by $f$. It is easy to
see that injectivity of $f$ implies that $V_a$ has no branch points. 

Let $W$ be a connected component of $V_a$. We claim that the image  
$f(W)$ is $(0,\infty)$.
It follows from Theorem~\ref{thm:milnor} and the above observation about
local maxima and minima that $f(W)$ is open. Now we show that $f$ is proper. 
Choose $w_n \in W$  with $f(w_n) \to b$ where $0 <b<
\infty$, so that by definition the lengths $ l_{\nu}(w_n)$ are bounded
above. Since $\t_{\mu} = a$ on $V_a$, it follows from Bridgeman's
inequality that the lengths $l_{\mu}(w_n)$ are also uniformly bounded
above. By Theorem~\ref{thm:limitpleating}, the corresponding sequence of
groups  $G(w_n)$ converge  algebraically to a group $G$ and since $G$ is clearly not Fuchsian, either
$G \in \P(\mu,\nu)$,    or $f(w_n)= l_{\nu}(w_n) \to 0$ or $l_{\mu}(w_n) \to 0$. By Proposition~\ref{prop:lengthbounds}, if $l_{\mu}(w_n) \to 0$ then $l_{\nu}(w_n) \to \infty$, ruling out the last possibility.
If $G \in \P(\mu,\nu)$ then by continuity $G$ corresponds to a point
$ w \in W$ with $f(w) = b$.  This shows that $f$ is proper and 
hence that $f(W) = (0,\infty)$ as claimed. 

Thus $W \cap L_c$ is non-empty for each $c \in
(0,\infty)$.  If $V_a$ had any other connected component, then for some $c$ the intersection $V_a \cap L_c$ would contain more than one point, which is impossible. The result follows.
\end{proof}

\section{Proof of the bending measure conjecture}
\label{sec:mainproof}

Suppose $a,b>0$. Any  quasifuchsian group with  $\beta^+ =
a\mu,\beta^- = b\nu $ necessarily lies in the variety $V_a$.
The following result therefore concludes the proof of  the bending measure conjecture:
 \begin{prop}
\label{prop:anglemono} Fix $a>0$. The angle $\t_{\nu}$ is strictly
monotonic on $V_a$.
\end{prop}
\begin{proof}
Let $q(c)$ denote the point $V_a \cap L_c$.
Choose sequences $\nu_n \to \nu, \mu_n \to \mu$, where as usual if $\mu \in \MLQ$ we assume that the sequence $\mu_n$
is constant. Let $q_n = q_n(c) \in
\P(\mu_n,\nu_n)$ be the unique point for which $l_{\nu_n} = c$ and 
$\t_{\mu_n} = a$; this exists for large $n$ since $\R(\mu_n,\nu_n)$ is close to $\R(\mu,\nu)$. By Bridgeman's inquality the lengths $l_{\nu_n}(q_n)$ are
uniformly bounded above. Hence by the limit pleating theorem, 
we can extract a subsequence for which $G(q_n)$ converges algebraically
to a group $G$.

We claim that  $l_{\mu_n} (q_n)$ does not tend to zero.
If $\mu \in \MLQ$ we can exactly follow exactly the argument in the second paragraph of the proof of Proposition~\ref{prop:lengthbounds}.
If $\mu  \notin \MLQ$, we proceed as in the third paragraph of that proof. Now $[\mu_n]$ is the projective class of the bending measure of $\bch^{+}(q_n)$ so that $l_{\mu_n^+}(q_n)$ and $l_{\mu_n^*}(q_n)$
coincide.  By Corollary~\ref{cor:thurstonbndry}, if $l_{\mu_n}(q_n) \to 0$  then 
the hyperbolic structures $\bch^{+}(q_n)$ converge to $[\mu]$.
This leads to a contradiction to the existence of the algebraic limit 
of the groups $G(q_n)$ exactly as in Proposition~\ref{prop:lengthbounds}.

We conclude from Theorem~\ref{thm:limitpleating} that $G$ is represented by a point $q' \in \P(\mu ,\nu)$. Moreover since $\t_{\nu_n}(q_n) \to \t_{\nu }(q'
)$ and $ l_{\nu_n}(q_n) \to l_{\nu }(q' )$, we deduce that $q' = q(c)$, so
that the limit is independent of the subsequence.

By definition $q_n(c)$ is the point $V^n_a \cap L^n_c$.
We already know that the angle function $\t_{\nu_n}(q_n(c))$
is monotonic in $c$. We have just shown that 
$q_n(c) \to q(c)$ and so $\t_{\nu_n}(q_n(c)) \to \t_{\nu}(q(c))$.
We deduce that $\t_{\nu}(q(c))$ is monotonic in $c$.

It remains to show that $\t_{\nu}(q(c))$ is strictly monotonic.
Now $\t_{\nu}(q(c))$ is a real analytic function of the 
parameters $l_{\mu}, l_{\nu}$ for $\P(\mu,\nu)$.  By
Theorem~\ref{thm:milnor}, on the real analytic variety $V_a$, each
length function is a real analytic function of some $t\in \RR$. Thus so
is $\t_{\nu}(q(c))$. We deduce that $\t_{\nu}(q(c))$ is either strictly
monotonic or constant on $V_a$.
To rule out the second possibility, notice that if  $\t_{\nu} $
is constant on $V_a$, then Bridgeman's inequality  
gives a uniform upper bound to $l_{\nu}$ on $V_a$, contradicting
Proposition~\ref{prop:shapeofV}.
\end{proof}

 \small{
}


\begin{thebibliography}{}

 \bibitem{Bon}
F.~Bonahon.
\newblock Kleinian groups which are almost Fuchsian.
 \newblock Preprint \texttt{arXiv:math:DG/0210233 v1}.
 
 \bibitem{BonS}
F.~Bonahon.
\newblock Shearing hyperbolic surfaces, bending pleated surfaces
and  Thurston's symplectic form.
\newblock {\em Ann. Fac. Sci. de Toulouse.}, 5(2):233--297, 
 1996. 

 
\bibitem{BO}
F.~Bonahon and J-P.~Otal.
\newblock Laminations  mesur\'ees de plissage des vari\'et\'es
hyperboliques
 de dimension 3.
 \newblock Annals of Mathematics, to appear. 
 
 
 \bibitem{Br}
M.~Bridgeman.
\newblock Average bending of convex pleated planes in $\HH^3$.
 \newblock {\em Inventiones} 132, 381--391, 1998. 

\bibitem{CEG}
 R.~D.~Canary, D.~B.~A.~Epstein and P.~Green.
\newblock Notes on notes of {T}hurston.
\newblock In D.~B.~A.~Epstein, editor, {\em Analytical and
Geometric Aspects of
  Hyperbolic Space}, LMS Lecture Notes 111, 3--92. Cambridge
University
  Press, 1987.

\bibitem{ChS}
Y.-E.~Choi and C.~Series. Lengths are coordinates for convex structures.
\newblock Warwick preprint, 2003.  


\bibitem{CEpG}
 R.~D.~Canary, D.~B.~A.~Epstein and P.~Green.
\newblock Notes on notes of {T}hurston.
\newblock In D.~B.~A.~Epstein, editor, {\em Analytical and
Geometric Aspects of
  Hyperbolic Space}, LMS Lecture Notes 111, 3--92. Cambridge
University
  Press, 1987.


\bibitem{EpM}
D.~B.~A. Epstein and A.~Marden.
\newblock Convex hulls in hyperbolic space, a theorem of
{S}ullivan, and
  measured pleated surfaces.
\newblock In D.~B.~A. Epstein, editor, {\em Analytical and
Geometric Aspects of
  Hyperbolic Space}, LMS Lecture Notes 111, 112--253. Cambridge
  University Press, 1987.

 
  
\bibitem{FLP}
A.~Fahti, P.~Laudenbach, and V.~Po{\'e}naru.
\newblock {\em Travaux de {T}hurston sur les surfaces}, 
  Ast{\'e}risque 66-67.
\newblock Soci{\'e}t{\'e} Math{\'e}matique de France, 1979.


 
\bibitem{KSconvex}
L.~Keen and C.~Series.
\newblock Continuity of convex hull boundaries.
\newblock {\em Pacific J. Math.}, 168(1):183--206, 1995.

\bibitem{KSbend}
L.~Keen and C.~Series.
\newblock How to bend pairs of punctured tori.
\newblock In J.~Dodziuk and L.~Keen, editors, {\em Lipa's
Legacy}, Contemp. Math. 211, 359--388.
  AMS, 1997.
  
  
  \bibitem{KSQF}
L.~Keen and C.~Series.
\newblock Pleating invariants for punctured torus groups.
\newblock {\em Topology} 43, 447--491, 2004.


\bibitem{KSRiley} L.~Keen and C.~Series.
\newblock The Riley slice of Schottky space,
\newblock {\em Proc. London Math. Soc.} 69,72--90, 1994.

 

\bibitem{KerckEA}
S.~Kerckhoff.
\newblock Earthquakes are analytic.
\newblock {\em Comment.\ Mat.\ Helv.}, 60:17--30, 1985.

\bibitem{KerckLM}
S.~Kerckhoff.
\newblock Lines of Minima in Teichm\"uller space.
\newblock {\em Duke Math J.}, 65:187--213, 1992. 

 \bibitem{KerckN}
S.~Kerckhoff.
\newblock The {N}ielsen realization problem.
\newblock {\em Ann. of Math}, 117(2):235--265, 1983.


\bibitem{Lec}
C.~Lecuire.
\newblock   Plissage des vari{\'e}t{\'e}s hyperboliques de dimension
$3$.  
  \newblock Preprint, 2002.


\bibitem{Mi}
J.~Milnor.
\newblock {\em Singular points of complex hypersurfaces.}
\newblock Annals of Math. Studies 61.
   Princeton University Press, 1968.



\bibitem{OtalH}
J.~P.~Otal.
\newblock {\em Le th\'eor{\`e}me d'hyperbolisation pour 
les vari{\'e}t{\'e}s fibr{\'e}es de dimension $3$}, 
  Ast{\'e}risque 235.
\newblock Soci{\'e}t{\'e} Math{\'e}matique de France, 1996.


 

\bibitem{PS}
J.~Parker and C.~Series.
 \newblock Bending formulae for convex hull boundaries.
 \newblock {\em J.d'Analyse Math.}, 67, 165--198, 1995.
 
\bibitem{SmB}
C.~Series.
\newblock	Limits of quasifuchsian groups
with small bending.
\newblock Duke J. Mathematics, to appear.  



 \bibitem{ThuN}
W.~Thurston.
\newblock {\em Geometry and Topology of Three-Manifolds}.
\newblock Princeton lecture notes, 1979.

 \bibitem{ThuII}
W.~Thurston.
\newblock Hyperbolic structures on $3$-manifolds {I}{I}.
\newblock  eprint at 
 {front.math.ucdavis.edu/search/author:Thurston+category:GT}.
 

\end{thebibliography}
\end{document}